\documentclass{amsart}
\newtheorem{thm}{Theorem}[section]
\newtheorem{cor}[thm]{Corollary}

\newtheorem{prop}[thm]{Proposition}
\theoremstyle{definition}
\newtheorem{defn}[thm]{Definition}
\theoremstyle{remark}
\newtheorem{rem}[thm]{Remark}
\numberwithin{equation}{section}

\begin{document}

\title[Classification of 4-dimensional nilpotent complex Leibniz
algebras.]{Classification of 4-dimensional nilpotent complex Leibniz algebras.}%

\author{S.Albeverio}
\address{[S.Albeverio] Institut f\"ur Angewandte Mathematik, Universit\"at
Bonn, Wegelestr.6, D-53115 Bonn (Germany).}

\author{B.A. Omirov}
\address{[B.A. Omirov] Institute of Mathematics of Academy of Uzbekistan. F.Khodjaev str, 29, 700125, Tashkent, Uzbekistan}
\email{omirovb@mail.ru}

\author{I.S.Rakhimov}
\address{[I.S.Rakhimov] University Putra Malaysia,
43400 UPM Serdang, Selangor (Malaysia)} \email{risamiddin@mail.ru,
isamiddin@science.upm.edu.my}



\begin{abstract}
The Leibniz algebras appeared as a generalization of the Lie
algebras. In this work we deal with the classification of
nilpotent complex Leibniz algebras of low dimensions. Namely, the
classification of nilpotent complex Leibniz algebras dimensions
less than 3 is extended to the dimension four.

\end{abstract}

\maketitle

\section{Introduction}

In \cite{L}, \cite{LF} several classes of new algebras were
introduced. Some of them have two generating operations and they
are called dialgebras. The first motivation to introduce such
algebraic structures was problems in algebraic K-theory and they
are related with well known Lie and associative algebras.

The categories of these algebras over their operads assemble into
the commutative diagram which reflects the Koszul duality of those
categories. The aim the present paper is to study structural
properties of one class of Loday's list, namely so called Leibniz
algebras.

Leibniz algebras present a "non commutative" (to be more precise,
a "non antisymmetric") analogue of Lie algebras and they were
introduced by J.-L.Loday \cite{L}, as algebras which satisfy the
following identity:

\begin{center}
$[x,[y,z]]=[[x,y],z]-[[x,z],y].$
\end{center}

They appeared to be related in a natural way to several topics
such as differential geometry, homological algebra, classic
algebraic topology, algebraic K-theory, loop spaces,
noncommutative geometry, etc.

Most papers concerning Leibniz algebras are devoted to the study
of homological problems \cite{CAS}, \cite{F}, \cite{P} and others.
The structure theory of Leibniz algebras mostly remains unexplored
and the extension of notions like simple and semisimple or radical
to Leibniz algebras have not been discussed.

Some structural results concerning nilpotency, classification of
low dimensional Leibniz algebras and related problems were
considered in \cite{AAO1},  \cite{AOR}, \cite{AO1}, \cite{AO2}.
For similar structural studies of Lie algebras see \cite{GK}.

The classification up to isomorphism of any class of algebras is a
fundamental and very difficult problem. It is one of the first
problems that one encounters when trying to understanding the
structure of a member of this class of algebras. In this paper we
are going to present the classification of complex nilpotent
Leibniz algebras of dimension four. From the geometrical point of
view the classification of a class of algebras corresponds to a
fibration of this class, the fiber being the isomorphic classes.
We will give representatives of each isomorphism class for
4-dimensional nilpotent complex Leibniz algebras. By the
classification here we mean the algebraic classification, i.e. the
determination of the types of isomorphic algebras, whereas
geometric classification is the problem of finding generic
structural constants in the sense of algebraic geometry. But the
geometric classification presupposes the algebraic classification.
Actually it is a prelude of our paper \cite{AOR}, where the
geometric classification problems of low dimensional complex
nilpotent Leibniz algebras were discussed.

We restrict our discussion to nilpotent Leibniz algebras of
dimension four since all Leibniz algebras of dimension less than
four already have been classified in \cite{AO1}, \cite{GK},
\cite{L}.

\begin{defn} An algebra L over a field F is said to be a
Leibniz algebra if it satisfies the following Leibniz identity:
\begin{center}
$[x,[y,z]]=[[x,y],z]-[[x,z],y],$
\end{center}
where [ , ] denotes the multiplication in L.
\end{defn}
\begin{rem} If a Leibniz algebra has the additional property
of antisymmetricity  $[x,x]=0$ then the Leibniz identity can be
easily reduced to the Jacobi identity:
$$[[x,y],z]+[[y,z],x]+[[z,x],y]=0.$$
\end{rem}

Therefore a Leibniz algebra is a generalization of well known Lie
algebras.

Another generalization of Lie algebras was given by A.Mal'tzev
\cite{S}. They are called Mal'tzev algebras and satisfy the
following two identities:

$$[x,y]=-[y,x],$$
$$[[[x,y],z],x]+[[[y,z],x],x]+[[[z,x],x],y]=[[x,y],[x,z]].$$

It is clear that the intersection of the varieties of Leibniz
algebras and Mal'tzev algebras coincides exactly with the variety
of Lie algebras.

All algebras considered in this paper suppose to be defined over
the field of complex numbers.

\section {Nilpotent Leibniz algebras.}

In this section we resume some results on nilpotent Leibniz
algebras. The concept of nilpotency in the case of Leibniz
algebras can be defined the similar manner like in Lie algebras
case.

Let $L$ be a complex Leibniz algebra. We put:
\begin{center}
$$L^1=L, L^{k+1}=[L^{k},L], \quad k \in N.$$
\end{center}

\begin{defn} A Leibniz algebra $L$ is said to be
nilpotent if there exists a natural $s \in N$, such that $L^s=0.$
\end{defn}

\begin{defn} An $n$-dimensional Leibniz algebra $L$ is
said to be nulfiliform if $dim L^i=n-i+1$, where $2 \leq i \leq
n+1.$
\end{defn}

\begin{defn} An $n$-dimensional Leibniz algebra $L$ is
said to be filiform if $dim L^i=n-i$, where $2 \leq i \leq n$.
\end{defn}

For a given $n-$dimensional nilpotent Leibniz algebra $L$ we
define the following isomorphism invariant:
\begin{center}
$\chi(L)=(dimL^1,dimL^2,\cdot\cdot\cdot,dimL^{n-1},dimL^n)$
\end{center}

It is evident that
\begin{center}
$dimL^1>dimL^2> \cdot\cdot\cdot > dimL^k > \cdot\cdot\cdot$
\end{center}

The following result can be easily proved.
\begin{prop} If for an n-dimensional Leibniz
algebra $L$ the first two coordinates of the invariant $\chi(L)$
are equal to n and n-1, respectively, then $L$ is a nulfiliform
Leibniz algebra.
\end{prop}

 Let $L$ be an associative algebra over the complex
numbers $C$ and let $A$ be the associative algebra obtained from
$L$ by the external joining of a unit, i.e. $A=L \oplus C$.

\begin{prop} Let $L$ be a finite dimensional
nilpotent associative algebra. Then the algebra $A=L \oplus C$
does not contain any nontrivial idempotents.
\end{prop}

\begin{proof} Let $a\neq0$ be an element of $L$ and $1+a$ be
an idempotent of $A=L+C$.  Let $m$ be the index of nilpotence of
$a$. It is evident that $1+a=(1+a)^n$ for any $ n \in N$. But
$(1+a)^n=1+ \sum\limits_{k=1}^n C_n^k a^k$ implies $a=na+
\sum\limits_{k=2}^n C_n^ka^k,$ and multiplying both sides of this
equality by $a^{m-2}$  we obtain $(n-1)a=0$, that shows $a=0.$
\end{proof}

\begin{cor} Let $L$ be a finite dimensional Leibniz
algebra with $L^3=0.$ Then the algebra $A=L \oplus C$ has no
nontrivial idempotents.
\end{cor}

\begin{prop} Let $L_1$ and $L_2$ be finite
dimensional associative algebras without unit. Then $L_1 \oplus C
\cong L_2 \oplus C$ if only if $L_1 \cong L_2$.\end{prop}

\begin{prop} \cite{AO2}. Up to isomorphism, there is only
one $n$-dimensional non Lie nulfiliform Leibniz algebra. It can be
given by the following table of multiplications:
$$[e_i,e_1]=e_{i+1}, \quad 1 \leq i \leq {n-1}$$
where $\{{e_1,e_2, \ldots,e_n}\}$ is a basis of $L$ and omitted
products are equal to zero.
\end{prop}
\begin{prop} \cite{AO2}. Any $(n+1)$-dimensional complex
non Lie filiform Leibniz algebra is isomorphic to one of the
following two non Lie filiform Leibniz algebras:
$$\mu_{1}^{\overline{\alpha},\theta}: [e_0,e_0]=e_2, \quad
[e_i,e_0]=e_{i+1} \quad (1 \leq i \leq n-1),$$
$$[e_0,e_1]=\alpha_3e_3 + \alpha_4e_4 + \cdot\cdot\cdot +
\alpha_{n-1} e_{n-1} + \theta e_n,$$
$$[e_i,e_1]=\alpha_3 e_{i+2} +\alpha_4 e_{i+3}+\cdot\cdot\cdot+\alpha_{n+1-i}e_n
 \quad (1\leq i \leq n-2)$$ \\
 and
$$\mu_{2}^{\overline{\beta},\gamma}: [e_0,e_0]=e_2, \quad [e_i,e_0]=e_{i+1} \quad (2 \leq i \leq n-1),$$
$$[e_0,e_1]=\beta_3 e_3+ \beta_4 e_4+\cdot\cdot\cdot+ \beta_n e_n, \quad
[e_1,e_1]=\gamma e_n,$$
$$[e_i,e_1]=\beta_3e_{i+2}+\beta_4e_{i+3}+\cdot\cdot\cdot+\beta_{n+1-i}e_n \quad (2 \leq i
\leq{n-2}),$$ where $\{{e_0,e_1,\cdot\cdot\cdot,e_n}\}$ is a basis
of $L$ and omitted products are equal to zero. \end{prop}

Note that the algebras $\mu_{1}^{\overline{\alpha},\theta}$ and
$\mu_{2}^{\overline{\beta},\gamma}$ are not isomorphic for any
parameters $\alpha_i, \theta, \beta_j, \gamma.$

\section{Description of 4-dimensional nilpotent complex
Leibniz algebras.}

In this section we are going to classify four-dimensional
nilpotent complex Leibniz algebras. Here we shall use the
classification of unital associative algebras of dimension five
\cite{M}. Namely, from that classification by some restrictions we
will select some Leibniz algebras in our list. Moreover, the
description of nul and filiform Leibniz algebras \cite{AO2} is
used as well.

\begin{prop} Any 4-dimensional nilpotent complex
Leibniz algebra $L$ belongs to one of the following types of
algebras:

(i)   nulfiliform Leibniz algebras, that is  $\chi(L)=(4,3,2,1)$

(ii)  filiform Leibniz algebras, that is $\chi(L)=(4,2,1,0)$

(iii) associative algebras, with  $\chi(L):=(4,2,0,0)$ or
$(4,1,0,0)$.

(iv)  abelian, that is $\chi(L)=(4,0,0,0)$.
\end{prop}

\begin{proof}
Let us consider all possible cases for $\chi(L)$:

If $\chi(L)=(4,3,2,1),$ then in this case $L$ is nulfiliform (see
proposition 2.8).

If $\chi(L)=(4,2,1,0),$ then in this case $L$ is filiform (see
proposition 2.9).

It is clear that in the cases when $\chi(L)$ is equal to (4,2,0,0)
and (4,1,0,0) respectively, then the algebra $L$ is associative.

In case $\chi(L)=(4,0,0,0)$ we have abelian algebra.

It is easy to see that proposition 2.4 implies that the cases
$\chi(L):=(4,3,2,0),$ $(4,3,1,0)$ and $(4,3,0,0)$ are not
possible. \end{proof}

Henceforth as a matter of convenience we assume that the undefined
multiplications are zero, we do not consider neither abelian
algebras, nor Lie algebras and split Leibniz algebras i.e. we do
not consider Leibniz algebras which are direct sums of proper
ideals.

Let $\{e_1,e_2,e_3,e_4\}$ be a basis of $L$.

\begin{thm} Up to isomorphism, there exist five one
parametric families and twelve concrete representatives of
nilpotent Leibniz algebras of dimension four, namely:
$$\Re_1: [e_1,e_1]=e_2, [e_2,e_1]=e_3,
[e_3,e_1]=e_4;$$
$$\Re_{2}: [e_1,e_1]=e_3, [e_1,e_2]=e_4,
[e_2,e_1]=e_3, [e_3,e_1]=e_4;$$
$$\Re_{3}: [e_1,e_1]=e_3, [e_2,e_1]=e_3, [e_3,e_1]=e_4;$$
$$\Re_4(\alpha): [e_1,e_1]=e_3, [e_1,e_2]=\alpha e_4,
[e_2,e_1]=e_3, [e_2,e_2]=e_4, [e_3,e_1]=e_4, \quad \alpha\in
\{0,1\};$$
$$\Re_{5}:  [e_1,e_1]=e_3, [e_1,e_2]=e_4, [e_3,e_1]=e_4;$$
$$\Re_6: [e_1,e_1]=e_3, [e_2,e_2]=e_4, [e_3,e_1]=e_4 ;$$
$$\Re_7: [e_1,e_1]=e_4, [e_1,e_2]=e_3, [e_2,e_1]=- e_3, [e_2,e_2]= -2e_3 +e_4$$
$$\Re_8: [e_1,e_2]=e_3, [e_2,e_1]=e_4, [e_2,e_2]=-e_3;$$
$$\Re_9(\alpha):  [e_1,e_1]=e_3, [e_1,e_2]=e_4, [e_2,e_1]=-\alpha e_3,
[e_2,e_2]=-e_4, \quad \alpha\in C;$$
$$\Re_{10}(\alpha):  [e_1,e_1]=e_4, [e_1,e_2]=\alpha e_4,
[e_2,e_1]=-\alpha e_4, [e_2,e_2]=e_4, [e_3,e_3]=e_4, \quad
\alpha\in C;$$
$$\Re_{11}: [e_1,e_2]=e_4, [e_1,e_3]=e_4, [e_2,e_1]=-e_4, [e_2,e_2]=e_4, [e_3,e_1]=e_4;$$
$$\Re_{12}: [e_1,e_1]=e_4, [e_1,e_2]=e_4, [e_2,e_1]=-e_4, [e_3,e_3]=e_4;$$
$$\Re_{13}: [e_1,e_2]=e_3, [e_2,e_1]=e_4;$$
$$\Re_{14}: [e_1,e_2]=e_3, [e_2,e_1]=-e_3, [e_2,e_2]=e_4;$$
$$\Re_{15}: [e_2,e_1]=e_4, [e_2,e_2]=e_3;$$
$$\Re_{16}(\alpha): [e_1,e_2]=e_4, [e_2,e_1]=(1+\alpha)/(1-\alpha)e_4,
[e_2,e_2]=e_3, \quad \alpha\in {C\backslash\{1\}};$$
$$\Re_{17}: [e_1,e_2]=e_4, [e_2,e_1]=-e_4, [e_3,e_3]= e_4.$$
\end{thm}
\begin{proof} By proposition 3.1 any 4-dimensional nilpotent
complex Leibniz algebra is either nulfiliform or filiform or
associative. Let us consider each of these classes separately.

Let $L$ be a nulfiliform. Then by proposition 2.8 there is only
one nulfiliform Leibniz algebra and it can be given by the table:
$$[e_1,e_1]=e_2, [e_2,e_1]=e_3, [e_3,e_1]=e_4.$$ Thus, in this case
we have the algebra $\Re_1$ in our list.

Next let $L$ be a filiform. Then proposition 2.9. implies that
there are two classes of filiform Leibniz algebras. In our case
they are:  $$\nabla(\alpha,\beta): [e_1,e_1]=e_3, [e_1,e_2]=\alpha
e_4, [e_2,e_1]=e_3, [e_2,e_2]=\beta e_4, [e_3,e_1]=e_4;$$ and
$$\Omega(\alpha,\beta): [e_1,e_1]=e_3, [e_1,e_2]=\alpha e_4,
[e_2,e_2]=\beta e_4, [e_3,e_1]=e_4.$$

Let us consider the case of the algebras $\nabla(\alpha,\beta).$

{\bf Case 1.1}: $\beta = 0$ and $\alpha \neq 0$. Then the
following change of the basis $\{e_1, e_2, e_3, e_4\}$ reduces
$\nabla(\alpha,\beta)$ to $\nabla(1,0)$:
$$e_1^{'}= \alpha e_1, e_2^{'}= \alpha e_2, e_3^{'}= \alpha^2
e_3, e_4^{'}= \alpha^3 e_4.$$

Thus $\nabla(1,0)$ coincides with $\Re_2$ in our list.

{\bf Case 1.2:} $\beta = 0$ and $\alpha = 0$. Then $\nabla(0,0)$
is the $\Re_3$.

{\bf Case 2:} $\beta \neq 0$. In this case by changing the basis
$\{e_1, e_2, e_3, e_4\}$ in the following way:
$$e_1^{'} = \beta e_1, e_2^{'} = \beta e_2, e_3^{'}= \beta^2 e_3, e_4^{'}= \beta^3 e_4$$
the above algebra $\nabla(\alpha,\beta)$ can be reduced to the
type $\nabla(\alpha,1)$.

Thus we have the following table for the algebras
$\nabla(\alpha,1)$: $$[e_1,e_1]=e_3, [e_1,e_2]=\alpha e_4,
[e_2,e_1]=e_3, [e_2,e_2]=e_4, [e_3,e_1]=e_4.$$ Let now pick out
the isomorphic classes within this family of algebras. It is easy
to see if we take the following possible changing of the basis
$\{e_{1},e_{2},e_{3},e_{4}\}$:

$$e_1^{'} = ae_1+be_{2}, e_2^{'}=(a+b)e_2+b(\alpha-1)e_{3},$$
$$e_3^{'}= a(a+b)e_3+ b(a\alpha+b)e_{4}, e_4^{'}= a^2(a+b)e_4$$
then for the parameters $\alpha$, $\alpha^{'},$ $a$, $b$ we get
the following equalities $\alpha^{'}=\frac{a\alpha+b}{a^{2}}$ and
$b=a^{2}-a$. It follows that at $\alpha=1$ we obtain
$\alpha^{'}=1$. But at $\alpha\neq1$ putting $a=1-\alpha$ we get
$\alpha^{'}=0$. Therefore up to isomorphism there are only two
algebras:
$$\Re_4(1): [e_1,e_1]=e_3, [e_2,e_1]=e_3, [e_2,e_2]=e_4, [e_3,e_1]=e_4, [e_1,e_2]=e_4;$$
$$\Re_4(0): [e_1,e_1]=e_3, [e_2,e_1]=e_3, [e_2,e_2]=e_4, [e_3,e_1]=e_4.$$
We note that since the dimensions of maximal abelian subalgebras
of the algebras $\Re_2$ and $\Re_4(\alpha)$ are different and
therefore $\Re_2$ never is isomorphic to the algebra
$\Re_4(\alpha)$ for any $\alpha\in\{0,1\}$. The algebra $\Re_3$ is
not isomorphic to the algebras $\Re_2$ and $\Re_4(\alpha)$ for any
$\alpha\in\{0,1\}$, because of the difference in dimensions of the
right annihilators.

Now let us consider the class $\Omega(\alpha,\beta)$.  There are
two possible cases.

{\bf Case 1:} $\beta = 0$ and $\alpha \neq 0$. Then the
transformation
$$e_1^{'}= e_1, e_2^{'}= \alpha^{-1}e_2, e_3^{'}=e_3, e_4^{'}=e_4$$
leads to $\Omega(1,0)$ which coincides with $\Re_5$.

{\bf Case 2}: $\beta \neq 0$. It is easy to see in this case the
transformation

$$e_1^{'} = \beta e_1, e_2^{'} =\beta e_2, e_3^{'}= \beta^2 e_3, e_4^{'}=\beta^3 e_4$$
leads to the type $$\Omega(\alpha,1): [e_1,e_1]=e_3,
[e_1,e_2]=\alpha e_4, [e_2,e_2]=e_4, [e_3,e_1]=e_4.$$ We consider
the following possible changing for basis within this class:

$$e_1^{'}=a e_1 + b e_2, \quad e_2^{'}=c e_2 - a{bc}^{-1}e_{3},$$
$$e_3^{'}=a^2 e_3 + b(a\alpha+b)e_{4}, \quad e_4^{'}=a^3e_4.$$

It is easy to see again after some calculations that for different
parameters $\alpha$ and $\alpha{'}$ we will get the equalities: $
\alpha{'}=\frac{c(a\alpha +b)}{a^3}$ and $c^{2}=a^{3}$ therefore
taking $b=-a\alpha$ we have that $\alpha=0$. Thus we obtain that
the algebras $\Omega(\alpha,1)$ are isomorphic to the algebra
$\Re_6$ with the following table of multiplications:
$$[e_1,e_1]=e_3, [e_2,e_2]=e_4, [e_3,e_1]=e_4.$$

It should be noted that the algebra $\Omega(0,0)$ is split.

It is easy to see that the algebra $\Re_5$ is not isomorphic to
the algebra $\Re_6,$ because of they have different dimensions of
left annihilators.

Now we will suppose that our algebra $L$ has the type iii) in
proposition 3.1, that is to say that $L$ is a Leibniz algebra with
$\chi(L):=$ (4,2,0,0) or (4,1,0,0), in particular $L$ is
associative. Then all results on associative algebras are
applicable to our case and henceforth we deal with associative
algebras under corresponding conditions on $L$.

We consider the associative algebra $A=L \oplus C$. It is
5-dimensional and unital. Due to G.Mazzola's paper \cite{M} we
have got the classification of all associative algebras with unit
in dimension five. Having the classification of G.Mazzola and
taking into a number of central idempotents it is easy to see that
our algebra $A$ is not isomorphic to the algebras with numbers
1-15, 25, 26, 38, 39, 55 in the list of G.Mazzola. Using corollary
2.6 and the fact that under isomorphism any idempotent is mapped
in to an idempotent it is easy to see that $A$ is not isomorphic
to the algebras 16-22, 27, 28, 29, 40, 41, 46, 47, 52, 58 in the
list of G.Mazzola. Moreover, the condition $L^3=0$ implies that
$A$ is not isomorphic to the algebras 23, 24, 33, 34, 44 in the
list of G.Mazzola. Then we will examine the others algebras in the
list of G.Mazzola and separate suitable cases that give the
required algebras.

N30. $A=C<x,y>/(xy+yx,xy-yx+y^2-x^2)+(x,y)^3$. Choose the basis of
$A$: $e_0=1, e_1=x, e_2=y, e_3=xy, e_4=x^2$. Then the subalgebra
$L=<e_1, e_2, e_3, e_4>$ with the following law gives the algebra:
$$\Re_7: [e_1,e_1]=e_4, [e_1,e_2]=e_3, [e_2,e_1]=-e_3, [e_2,e_2]=-2e_3+e_4.$$

N31. $A=C<x,y>/(x^2,xy+y^2)+(x,y)^3$. Choose the basis of $A$:
$e_0=1, e_1=x, e_2=y, e_3=xy, e_4=yx$. Then the subalgebra
$L=<e_1, e_2, e_3, e_4>$ with the following law gives the
corresponding algebra:
$$\Re_8: [e_1,e_2]=e_3, [e_2,e_1]=e_4, [e_2,e_2]=-e_3.$$

N32. $A=C<x,y>/(xy+y^2,\alpha x^2+yx)+(x,y)^3$. Choose the basis
of $A$:  $e_0=1, e_1=x, e_2=y, e_3=x^2, e_4=xy$. Then the
subalgebra $L=<e_1, e_2, e_3, e_4>$ with the following law gives
the algebra:
$$\Re_9(\alpha): [e_1,e_1]=e_3, [e_1,e_2]=e_4, [e_2,e_1]=-\alpha e_3, [e_2,e_2]=-e_4.$$
Note that by proposition 2.7 for  different $\alpha_1$ and
$\alpha_2$ the algebras $\Re_9(\alpha_1)$ and $\Re(\alpha_2)$ are
not isomorphic.

N35. $A=C<x,y,z>/(xz,yz,zx,zy,x^2-y^2,x^2-z^2,xy+yx,\alpha
x^2+yx)$, where $(\alpha\neq0)$. Choose the basis of $A$:  $e_0=1,
e_1=x, e_2=y, e_3=z, e_4=x^2$. Then the subalgebra $L=<e_1, e_2,
e_3, e_4>$ with the following law gives the corresponding algebra:
$$\Re_{10}(\alpha): [e_1,e_1]=e_4, [e_1,e_2]=\alpha e_4, [e_2,e_1]=-\alpha e_4, [e_2,e_2]=e_4,
[e_3,e_3]=e_4.$$

We remark that the algebras $\Re_{10}(\alpha_1)$ and
$\Re_{10}(\alpha_2)$ $(\alpha_1 \neq \alpha_2)$ are not isomorphic
except for the case where $\alpha_2=-\alpha_1$, in which case they
are indeed isomorphic \cite{M}.

N36. $A=C<x,y,z>/(x^2,yz,zy,z^2,xy-xz,xy+yx,yx+y^2,yx+zx)$. As a
basis here we take $e_0=1, e_1=x, e_2=y, e_3=z, e_4=xy$. Then the
subalgebra $L=<e_1, e_2, e_3, e_4>$ with the following law gives
the corresponding algebra:
$$\Re_{11}: [e_1,e_2]=e_4, [e_1,e_3]=e_4, [e_2,e_1]=-e_4, [e_2,e_2]=e_4, [e_3,e_1]=e_4.$$

N37. $A=C<x,y,z>/(xz, y^2,yz,zx,zy,x^2-z^2,x^2-xy,x^2+yx)$. As a
basis we take $e_0=1, e_1=x, e_2=y, e_3=z, e_4=xy$. Then the
subalgebra $L=<e_1, e_2, e_3, e_4>$ with the following law gives
the algebra:
$$\Re_{12}: [e_1,e_1]=e_4, [e_1,e_2]=e_4, [e_2,e_1]=-e_4, [e_3,e_3]=e_4.$$

N42. $A=C<x,y>/(x^2,y^2)+(x,y)^3$.  As a basis we take $e_0=1,
e_1=x, e_2=y, e_3=xy, e_4=yx$. Then the subalgebra $L=<e_1, e_2,
e_3, e_4>$  gives the algebra:
$$\Re_{13}: [e_1,e_2]=e_3, [e_2,e_1]=e_4.$$

N43. $A=C[x,y]/(y^3,xy,x^3)$.  As a basis we take $e_0=1, e_1=x,
e_2=y, e_3=x^2, e_4=y^2$. Then the subalgebra $L=<e_1, e_2, e_3,
e_4>$  gives the algebra:

$$[e_1,e_1]=e_3, [e_2,e_2]=e_4.$$
It is obviously that this algebra is decomposable.

N48. $A=C<x,y>/(x^2,xy+yx)+(x,y)^3$.  As a basis we take $e_0=1,
e_1=x, e_2=y, e_3=xy, e_4=y^2$. Then the subalgebra $L=<e_1, e_2,
e_3, e_4>$ coincides with the algebra:
$$\Re_{14}: [e_1,e_2]=e_3, [e_2,e_1]=-e_3, [e_2,e_2]=e_4.$$

N49($\alpha=1$). $A=C<x,y>/(x^2,xy)+(x,y)^3$.  As a basis we take
$e_0=1, e_1=x, e_2=y, e_3=y^2, e_4=yx$. Then the subalgebra
$L=<e_1, e_2, e_3, e_4>$  gives the algebra:
$$\Re_{15}: [e_2,e_1]=e_4, [e_2,e_2]=e_3.$$

N49($\alpha \neq 1$).
$A=C<x,y>/(x^2,(1+\alpha)xy+(1-\alpha)yx)+(x,y)^3$.  As a basis we
take $e_0=1, e_1=x, e_2=y, e_3=y^2, e_4=xy$. Then the subalgebra
$L=<e_1, e_2, e_3, e_4>$ with the following law gives the
corresponding algebra:
$$\Re_{16}(\alpha): [e_1,e_2]=e_4, [e_2,e_1]=\frac{(1+\alpha)}{(1-\alpha)}e_4, [e_2,e_2]=e_3.$$

Note that again by proposition 2.7 for different values of
$\alpha$ we obtain non isomorphic algebras.

N50. $A=C<x,y,z>/(x^2,xz,y^2,yz,zx,zy,xy+yx,yx+z^2)$.  As a basis
we take $e_0=1, e_1=x, e_2=y, e_3=z, e_4=xy$. Then the subalgebra
$L=<e_1, e_2, e_3, e_4>$ gives the algebra:
$$\Re_{17}: [e_1,e_2]=e_4, [e_2,e_1]=-e_4, [e_3,e_3]=e_4.$$

N51. $A=C<x,y,z>/(xz, yz, zx, zy, x^2-y^2, x^2-z^2, xy, yx)$.  As
a basis we take $e_0=1, e_1=x, e_2=y, e_3=z, e_4=x^2$. Then the
subalgebra $L=<e_1, e_2, e_3, e_4>$ has the following table:
$$[e_1,e_1]=e_4, [e_2,e_2]=e_4, [e_3,e_3]=e_4.$$
But this algebra can be included in to the list of algebras
$\Re_{10}(\alpha)$ for $\alpha = 0$.

Using proposition 2.7 it is easy to check that all obtained
algebras are pairwise not isomorphic.
\end{proof}

\begin{rem} All other algebras from the G.Mazzola's list
are either Lie algebras or split Leibniz algebras.
\end{rem}

It should be noted that unlike the case of Lie algebras, where in
dimension four there is only one non split algebra, we have got a
lot of Leibniz algebras in dimension four.

Summing the classification of the above theorem and the
classifications of complex nilpotent Lie algebras dimensions less
that five and complex nilpotent Leibniz algebras of dimensions
less than four, we obtain the complete classification of complex
nilpotent Leibniz algebras dimension less than five.

{\it {\bf Acknowledgment.} The authors would like to convey their
sincere thanks to prof. Sh.A. Ayupov for his comprehensive
assistance and support in the completion of this paper. The second
and third named authors would like to acknowledge the hospitality
to the "Institut f\"ur Angewandte Mathematik", Universit\"at Bonn
(Germany). This work is supported in part by the DFG 436 USB 113/4
project (Germany) and the Fundamental Science Foundation of
Uzbekistan.}


\begin{thebibliography}{12}

\bibitem{AAO1}
{\rm Albeverio S., Ayupov. Sh.A., Omirov B.A.} {\it On nilpotent
and simple Leibniz algebras.} Comm. in Algebra, v. 33(1), 2005, p.
159-172.

\bibitem{AOR}
{\rm Albeverio S., Omirov B.A., Rakhimov I.S.} {\it Varieties of
nilpotent complex Leibniz algebras of dimension less than five.}
Comm. in Algebra, v. 33(5), 2005, p. 1575-1585.

\bibitem{AO1}
{\rm Ayupov Sh.A., Omirov B.A.} {\it On Leibniz algebras.} Algebra
and operators theory, Proceeding of the Colloquium in Tashkent
1997. Kluwer Academic Publishers 1998, p. 1-13.

\bibitem{AO2}
{\rm Ayupov Sh.A., Omirov B.A.} {\it On some classes of nilpotent
Leibniz algebras.} Siberian Math. Journal, v.42(1), 2001, p.
18-29.

\bibitem{CAS}
{\rm Casas J.M., Pirashvili T.} {\it Ten-term exact sequence of
Leibniz homology.} J. Algebra, v. 231, 2000, p. 258-264.

\bibitem{F}
{\rm Frabetti A.} {\it Leibniz homology of dialgebras of
matrices,} Journ. Pure Appl. Alg., v. 129, 1998, p. 123-141.

\bibitem{GK}
{\rm Goze M., Khakimdjanov Yu.} {\it Nilpotent Lie Algebras.}
Kluwer Academic publishers, Dordrecht, 1996, 336 p.

\bibitem{L}
{\rm Loday J.-L.} {\it Une version non commutative des
alg$\grave{e}$bres de Lie: les alg$\grave{e}$bres de Leibniz,}
Ens. Math., v. 39, 1993, p. 269-293.

\bibitem{LF}
{\rm Loday J.-L., Frabetti A., Chapoton F., Goichot F.} {\it
Dialgebras and Related Operads,} Lect. Notes in Math., v. 1763
2001, 133 p.

\bibitem{M}
{\rm Mazzola G.} {\it The algebraic and geometric classification
of associative algebras of dimension five.} Manuscripta math. 27,
1979, p. 1-21.

\bibitem{P}
{\rm Pirashvili T.} {\it On Leibniz homology,} Ann.Inst. Fourier,
v. 44(2), 1994, p. 401-411.

\bibitem{S}
{\rm Sagle A.A.} {\it Mal'tzev algebras.} Trans.Amer.Math.Soc., v.
101(3), 1961, p. 426-458.

\end{thebibliography}
\end{document}